\documentclass{amsart}
\usepackage{graphicx}
\usepackage{latexsym}
\usepackage{amsfonts}
\usepackage[all]{xy}
\usepackage{amssymb, mathrsfs, amsfonts, amsmath}
\usepackage{amsbsy}
\usepackage{amsfonts}
\newtheorem*{acknowledgement}{Acknowledgement}
\newtheorem{corollary}{Corollary}

\newtheorem{theorem}{Theorem}

\numberwithin{equation}{section}

\begin{document}
\title[Einstein warped products]{On the structure of Einstein warped product semi-Riemannian manifolds}
\author{Benedito Leandro$^{1}$}
\address{$^{1}$ Universidade Federal de Goi\'as, Centro de Ci\^encias Exatas, Regional Jata\'i, BR 364, km 195, 3800, 75801-615, Jata\'i, GO, Brazil.}
 \email{bleandroneto@gmail.com$^{1}$}
\author{M\'arcio Lemes de Sousa $^{2}$}
\address{$^{2}$ ICET-CUA, Universidade Federal de Mato Grosso, Av. Universit\'aria number 3500,
Pontal do Araguaia, MT, Brazil.}
 \email{marciolemesew@yahoo.com.br}
\author{Romildo Pina  $^{3}$}
\address{$^{3}$ Universidade Federal de Goi\'as, IME, 131, 74001-970, Goi\^ania, GO, Brazil.}
 \email{romildo@ufg.br$^{3}$}
\date{August 14, 2017.}

\begin{abstract}
In this paper we consider a class of Einstein warped product semi-Riemannian manifolds $\widehat{M}=M^{n}\times_{f}N^{m}$ with $n\geq3$ and $m\geq2$. For $\widehat{M}$ with compact base and Ricci-flat fiber, we prove that $\widehat{M}$ is simply a Riemannian product space. Then, when the base $M$ is conformal to a pseudo-Euclidean space which is invariant under the action of a $(n-1)$-dimensional translation group, we classify all such spaces. Furthermore, we get new examples of complete Einstein warped products Riemannian manifolds.
\end{abstract}
\maketitle

\section{Introduction}

Einstein manifolds are related with many questions in geometry and physics, for instance: Riemannian functionals and their critical points, Yang-Mills theory, self-dual manifolds of dimension four, exact solutions for the Einstein equation field. Today we already have in our hands many examples of Einstein manifolds, even the Ricci-flat ones (see \cite{besse,Oneil,LeandroPina,Romildo}). However, finding new examples of Einstein metrics is not an easy task. A common tool to make new examples of Einstein spaces is to consider warped product metrics (see \cite{LeandroPina,Romildo}). 

In \cite{besse}, a question was made about Einstein warped products:
\begin{eqnarray}\label{question}
\mbox{``Does there exist a compact Einstein warped}\nonumber\\
\mbox{product with nonconstant warping
function?"}
\end{eqnarray}

Inspired by the problem (\ref{question}), several authors explored this subject in an attempt to get examples of such manifolds. Kim and Kim \cite{kimkim} considered a compact Riemannian Einstein warped product with nonpositive scalar curvature. They proved that a manifold like this is just a product manifold.  Moreover, in \cite{BRS,Case1}, they considered (\ref{question}) without the compactness assumption. Barros, Batista and Ribeiro Jr \cite{BarrosBatistaRibeiro} also studied (\ref{question}) when the Einstein product manifold is complete and noncompact with nonpositive scalar curvature. It is worth to say that Case, Shu and Wei \cite{Case} proved that a shrinking quasi-Einstein metric has positive scalar curvature. Further, Sousa and Pina \cite{Romildo} were able to classify some structures of Einstein warped product on semi-Riemannian manifolds, they considered, for instance, the case in which the base and the fiber are Ricci-flat semi-Riemannian manifolds. Furthermore, they provided a  classification for a noncompact Ricci-flat warped product semi-Riemannian manifold with $1$-dimensional fiber, however the base is not necessarily a Ricci-flat manifold. More recently, Leandro and Pina \cite{LeandroPina} classified the static solutions for the vacuum Einstein equation field with cosmological constant not necessarily identically zero, when the base is invariant under the action
of a translation group. In particular, they provided a necessarily condition for integrability of the system of differential equations given by the invariance of the base for the static metric. 

When the base of an Einstein warped product is a compact Riemannian manifold and the fiber is a Ricci-flat semi-Riemannian manifold, we get a partial answer for (\ref{question}). Furthermore, when the base is not compact, we obtain new examples of Einstein warped products.
 
Now, we state our main results. 

\begin{theorem}\label{teo1}
Let $(\widehat{M}^{n+m}, \hat{g})=(M^{n},g)\times_{f}(N^{m},\tilde{g})$, $n\geq3$ and $m\geq2$, be an Einstein warped product semi-Riemannian manifold (non Ricci-flat), where $M$ is a compact Riemannian manifold and $N$ is a Ricci-flat semi-Riemannian manifold. Then $\widehat{M}$ is a product manifold, i.e., $f$ is trivial.
\end{theorem}

It is very natural to consider the next case (see Section \ref{SB}).

\begin{theorem}\label{teo2}
Let $(\widehat{M}^{n+m}, \hat{g})=(M^{n},g)\times_{f}(N^{m},\tilde{g})$, $n\geq3$ and $m\geq2$, be an Einstein warped product semi-Riemannian manifold (i.e., $\widehat{R}ic=\lambda\hat{g}$; $\lambda\neq0$), where $M$ is a compact Riemannian manifold with scalar curvature $R\leq\lambda(n-m)$, and $N$ is a semi-Riemannian manifold. Then $\widehat{M}$ is a product manifold, i.e., $f$ is trivial. Moreover, if the equality holds, then $N$ is Ricci-flat.
\end{theorem}

Now, we consider that the base is a noncompact Riemannian manifold. The next result was inspired, mainly, by Theorem \ref{teo2} and \cite{LeandroPina}, and gives the relationship between Ricci tensor $\widehat{R}ic$ of the warped metric $\hat{g}$ and the Ricci tensor $Ric$ for the metric of the base $g$.

\begin{theorem}\label{teo3b}
Let $(\widehat{M}^{n+m}, \hat{g})=(M^{n},g)\times_{f}(N^{m},\tilde{g})$, $n\geq3$ and $m\geq2$, be an Einstein warped product semi-Riemannian manifold (i.e., $\widehat{R}ic=\lambda\hat{g}$), where $M$ is a noncompact Riemannian manifold with constant scalar curvature $\lambda=\frac{R}{n-1}$, and $N$ is a semi-Riemannian manifold. Then $M$ is Ricci-flat if and only if the scalar curvature $R$ is zero.
\end{theorem}

Considering a conformal structure for the base of an Einstein warped product semi-Riemannian manifold, we have the next results. Furthermore, the following theorem is very technical. We consider that the base for such Einstein warped product manifold is conformal to a pseudo-Euclidean space which is invariant under the action of a $(n-1)$-dimensional translation group, and that the fiber is a Ricci-flat space. In order, for the reader to have a more intimate view of the next results, we recommend a previous reading of Section \ref{CFSI}.

\begin{theorem}\label{teo3a}
Let $(\widehat{M}^{n+m}, \hat{g})=(\mathbb{R}^{n}, \bar{g})\times_{f}(N^{m},\tilde{g})$, $n\geq3$ and $m\geq2$, be an warped product semi-Riemannian manifold such that $N$ is a Ricci-flat semi-Riemannian manifold. Let $(\mathbb{R}^{n}, g)$ be a pseudo-Euclidean space with coordinates $x =(x_{1},\ldots , x_{n})$ and $g_{ij} = \delta_{ij}\varepsilon_{i}$, $1\leq i,j\leq n$, where $\delta_{ij}$ is the delta Kronecker and $\varepsilon_{i} = \pm1$, with at least one $\varepsilon_{i} = 1$. Consider smooth functions $\varphi(\xi)$ and $f(\xi)$, where $\xi=\displaystyle\sum_{k=1}^{n}\alpha_{k}x_{k}$, $\alpha_{k}\in\mathbb{R}$, and $\displaystyle\sum_{k=1}^{n}\varepsilon_{k}\alpha_{k}^{2}=\kappa$. Then $(\mathbb{R}^{n}, \bar{g})\times_{f}(N^{m},\tilde{g})$, where $\bar{g}=\frac{1}{\varphi^{2}}g$, is an Einstein warped product semi-Riemannain manifold (i.e., $\widehat{R}ic=\lambda\hat{g}$) such that $f$ and $\varphi$ are given by:
\begin{eqnarray}\label{system2}
\left\{
\begin{array}{lcc}
(n-2)\varphi\varphi''-m\left(G\varphi\right)'=mG^{2}\\\\
\varphi\varphi''-(n-1)(\varphi')^{2}+mG\varphi'=\kappa\lambda \\\\
nG\varphi'-(G\varphi)'-mG^{2}=\kappa\lambda,
\end{array}
\right.
\end{eqnarray}
and
\begin{eqnarray}\label{sera3}
f=\Theta\exp\left(\int\frac{G}{\varphi}d\xi\right),
\end{eqnarray}
where $\Theta\in\mathbb{R}_{+}\backslash\{0\}$, $G(\xi)=\pm\sqrt{\frac{\kappa[\lambda(n-m)-\bar{R}]}{m(m-1)}}$ and $\kappa=\pm1$. Here $\bar{R}$ is the scalar curvature for $\bar{g}$.
\end{theorem}

The next result is a consequence of Theorem \ref{teo3a}.

\begin{theorem}\label{teo3}
Let $(\widehat{M}^{n+m}, \hat{g})=(M^{n},g)\times_{f}(N^{m},\tilde{g})$, $n\geq3$ and $m\geq2$, be an Einstein warped product semi-Riemannian manifold, where $M$ is conformal to a pseudo-Euclidean space invariant under the action of a $(n-1)$-dimensional translation group with constant scalar curvature (possibly zero), and $N$ is a Ricci-flat semi-Riemannian manifold. Then, $\widehat{M}$ is either
\begin{enumerate}
\item[(1)] a Ricci-flat semi-Riemannain manifold $(\mathbb{R}^{n},g)\times_{f}(N^{m},\tilde{g})$, such that $(\mathbb{R}^{n},g)$ is the pseudo-Euclidean space with warped function $f(\xi)=\Theta\exp{(A\xi)}$, where $\Theta>0,$ $A\neq0$ are nonnull constants, or\\
\item[(2)] conformal to $(\mathbb{R}^{n},g)\times(N^{m},\tilde{g})$, where $(\mathbb{R}^{n},g)$ is the pseudo-Euclidean space. The conformal function $\varphi$ is given by
\begin{eqnarray*}
\varphi(\xi)= \frac{1}{(-G\xi+C)^{2}};\quad\mbox{where}\quad G\neq0, C\in\mathbb{R}.
\end{eqnarray*} 
\end{enumerate}
Moreover, the conformal function is defined for $\xi\neq\frac{C}{G}$.
\end{theorem}
It is worth mentioning that the first item of Theorem \ref{teo3} was not considered in \cite{Romildo}.

From Theorem \ref{teo3} we can construct examples of complete Einstein warped product Riemannian manifolds.

\begin{corollary}\label{coro1}
Let $(N^{m},\tilde{g})$ be a complete Ricci-flat Riemannian manifold and $f(\xi)=\Theta\exp{(A\xi)}$, where $\Theta>0$ and $A\neq0$ are constants. Therefore, $(\mathbb{R}^{n},g_{can})\times_{f}(N^{m},\tilde{g})$ is a complete Ricci-flat warped product Riemannian manifold.
\end{corollary}

\begin{corollary}\label{coro2}
Let $(N^{m},\tilde{g})$ be a complete Ricci-flat Riemannian manifold and $f(x)= \frac{1}{x_{n}}$ with $x_{n}>0$. Therefore, $(\widehat{M},\hat{g})=(\mathbb{H}^{n},g_{can})\times_{f}(N^{m},\tilde{g})$ is a complete Riemannian Einstein warped product such that $$\widehat{R}ic=-\frac{m+n-1}{n(n-1)}\hat{g}.$$
\end{corollary}

The paper is organized as follows. Section \ref{SB} is divided in two subsections, namely, {\it General formulas} and {\it A conformal structure for the warped product with Ricci-flat fiber}, where will be provided the preliminary results. Further, in Section \ref{provas}, we will prove our main results.

\section{Preliminar}\label{SB}

Consider $(M^{n}, g)$ and $(N^{m},\tilde{g})$, with $n\geq3$ and $m\geq2$, semi-Riemannian manifolds, and let $f:M^{n}\rightarrow(0,+\infty)$ be a smooth function, the warped product $(\widehat{M}^{n+m},\hat{g})=(M^{n},g)\times_{f}(N^{m},\tilde{g})$ is a product manifold $M\times N$ with metric 
\begin{eqnarray*}
\hat{g}=g+f^{2}\tilde{g}.
\end{eqnarray*}
From Corollary 43 in \cite{Oneil}, we have that (see also \cite{kimkim})
\begin{eqnarray}\label{test1}
\widehat{R}ic=\lambda\hat{g}\Longleftrightarrow\left\{
\begin{array}{lcc}
Ric-\frac{m}{f}\nabla^{2}f=\lambda g\\
\widetilde{R}ic=\mu\tilde{g}\\
f\Delta f+(m-1)|\nabla f|^{2}+\lambda f^{2}=\mu
\end{array}
,\right.
\end{eqnarray}
where $\lambda$ and $\mu$ are constants. Which means that $\widehat{M}$ is an Einstein warped product if and only if (\ref{test1}) is satisfied. Here $\widehat{R}ic$, $\widetilde{R}ic$ and $Ric$ are, respectively, the Ricci tensor for $\hat{g}$, $\tilde{g}$ and $g$. Moreover, $\nabla^{2}f$, $\Delta f$ and $\nabla f$ are, respectively, the Hessian, The Laplacian and the gradient of $f$ for $g$.

\subsection{General formulas}\label{GF}

We derive some useful formulae from system (\ref{test1}). Contracting the first equation of (\ref{test1}) we get
\begin{eqnarray}\label{01}
Rf^{2}-mf\Delta f=nf^{2}\lambda,
\end{eqnarray}
where $R$ is the scalar curvature for $g$. From the third equation in (\ref{test1}) we have
\begin{eqnarray}\label{02}
mf\Delta f+m(m-1)|\nabla f|^{2}+m\lambda f^{2}=m\mu.
\end{eqnarray}
Then, from (\ref{01}) and (\ref{02}) we obtain
\begin{eqnarray}\label{oi}
|\nabla f|^{2}+\left[\frac{\lambda(m-n)+R}{m(m-1)}\right]f^{2}=\frac{\mu}{(m-1)}.
\end{eqnarray}

When the base is a Riemannian  manifold and the fiber is a Ricci-flat semi-Riemannian manifold (i.e., $\mu=0$), from (\ref{oi}) we obtain
\begin{eqnarray}\label{eqtop}
|\nabla f|^{2}+\left[\frac{\lambda(m-n)+R}{m(m-1)}\right]f^{2}=0.
\end{eqnarray}
Then, either
\begin{eqnarray*}
R\leq\lambda(n-m)
\end{eqnarray*}
or $f$ is trivial , i.e., $\widehat{M}$ is a product manifold.

\subsection{A conformal structure for the Warped product with Ricci-flat fiber}\label{CFSI}

In what follows, consider $(\mathbb{R}^{n}, g)$ and $(N^{m},\tilde{g})$ semi-Riemannian manifolds, and let $f:\mathbb{R}^{n}\rightarrow(0,+\infty)$ be a smooth function, the warped product $(\widehat{M}^{n+m},\hat{g})=(\mathbb{R}^{n},g)\times_{f}(N^{m},\tilde{g})$ is a product manifold $\mathbb{R}^{n}\times N$ with metric 
\begin{eqnarray*}
\hat{g}=g+f^{2}\tilde{g}.
\end{eqnarray*}

Let $(\mathbb{R}^{n}, g)$, $n\geq3$, be the standard pseudo-Euclidean space with metric $g$ and coordinates $(x_{1},\ldots,x_{n})$ with $g_{ij}=\delta_{ij}\varepsilon_{i}$, $1\leq i,j\leq n$, where $\delta_{ij}$ is the delta Kronecker, $\varepsilon_{i} = \pm1$, with at least one $\varepsilon_{i} = 1$. Consider $(\widehat{M}^{n+m},\hat{g})=(\mathbb{R}^{n}, \bar{g})\times_{f}(N^{m},\tilde{g})$ a warped product, where $\varphi:\mathbb{R}^{n}\rightarrow\mathbb{R}\backslash\{0\}$ is a smooth function such that $\bar{g}=\frac{g}{\varphi^{2}}$. Furthermore, we consider that $\widehat{M}$ is an Einstein semi-Riemannian manifold, i.e., $$\widehat{R}ic=\lambda\hat{g},$$
where $\widehat{R}ic$ is the Ricci tensor for the metric $\hat{g}$ and $\lambda\in\mathbb{R}$.

We use invariants for the group action (or subgroup) to reduce a partial differential equation into a system of ordinary differential equations \cite{olver}. To be more clear, we consider that $(\widehat{M}^{n+m},\hat{g})=(\mathbb{R}^n,\bar{g})\times_{f}(N^{m},\tilde{g})$ is such that the base is invariant under the action of a $(n-1)$-dimensional translation group (\cite{BarbosaPinaKeti,olver,LeandroPina,Romildo,Tenenblat}). More precisely, let $(\mathbb{R}^{n}, g)$ be the standard pseudo-euclidean space with metric $g$ and coordinates $(x_{1}, \cdots, x_{n})$, with $g_{ij} = \delta_{ij}\varepsilon_{i}$, $1\leq i, j\leq n$, where $\delta_{ij}$ is the delta Kronecker,
$\varepsilon_{i} = \pm1$, with at least one $\varepsilon_{i} = 1$. Let $\xi=\displaystyle\sum_{i}\alpha_{i}x_{i}$, $\alpha_{i}\in\mathbb{R}$, be a basic invariant for a $(n-1)$-dimensional translation group where $\alpha=\displaystyle\sum_{i}\alpha_{i}\frac{\partial}{\partial x_{i}}$ is a timelike, lightlike or spacelike vector, i.e., $\displaystyle\sum_{i}\varepsilon_{i}\alpha_{i}^{2}=-1,0,$ or $1$, respectively. Then we consider $\varphi(\xi)$ and $f(\xi)$ non-trivial differentiable functions such that 
\begin{eqnarray*}
\varphi_{x_{i}}=\varphi'\alpha_{i}\quad\mbox{and}\quad f_{x_{i}}=f'\alpha_{i}.
\end{eqnarray*}

Moreover, it is well known (see \cite{BarbosaPinaKeti,LeandroPina,Romildo}) that if $\bar{g}=\frac{1}{\varphi^{2}}g$, then the Ricci tensor $\bar{R}ic$ for $\bar{g}$ is given by
$$\bar{R}ic=\frac{1}{\varphi^{2}}\{(n-2)\varphi\nabla^{2}\varphi + [\varphi\Delta\varphi - (n-1)|\nabla\varphi|^{2}]g\},$$ where $\nabla^{2}\varphi$, $\Delta\varphi$ and $\nabla\varphi$ are, respectively, the Hessian, the Laplacian and the gradient of $\varphi$ for the metric $g$.
Hence, the scalar curvature of $\bar{g}$ is given by
\begin{eqnarray}\label{scalarcurvature}
\bar{R}&=&\displaystyle\sum_{k=1}^{n}\varepsilon_{k}\varphi^{2}\left(\bar{R}ic\right)_{kk}=(n-1)(2\varphi\Delta\varphi - n|\nabla\varphi|^{2})\nonumber\\
&=&(n-1)[2\varphi\varphi''-n(\varphi)^{2}]\displaystyle\sum_{i}\varepsilon_{i}\alpha_{i}^{2}.
\end{eqnarray}
In what follows, we denote $\kappa=\displaystyle\sum_{i}\varepsilon_{i}\alpha_{i}^{2}$.

When the fiber $N$ is a Ricci-flat semi-Riemannian manifold, we already know from Theorem 1.2 in \cite{Romildo} that $\varphi(\xi)$ and $f(\xi)$ satisfy the following system of differential equations 

\begin{eqnarray}\label{system}
\left\{
\begin{array}{lcc}
(n-2)f\varphi''-mf''\varphi-2m\varphi'f'=0;\\\\
f\varphi\varphi''-(n-1)f(\varphi')^{2}+m\varphi\varphi'f'=\kappa\lambda f;\\\\
(n-2)f\varphi\varphi'f'-(m-1)\varphi^{2}(f')^{2}-ff''\varphi^{2}=\kappa\lambda f^{2}.
\end{array}
\right.
\end{eqnarray}
Note that the case where $\kappa=0$ was proved in \cite{Romildo}. Therefore, we only consider the case $\kappa=\pm1$.

\

\section{Proof of the main results}\label{provas}

\

\noindent {\bf Proof of Theorem \ref{teo1}:}
In fact, from the third equation of the system (\ref{test1}) we get that
\begin{eqnarray}\label{kimkimeq}
div\left(f\nabla f\right)+(m-2)|\nabla f|^{2}+\lambda f^{2}=\mu.
\end{eqnarray}
Moreover, if $N$ is Ricci-flat, from (\ref{kimkimeq}) we obtain
\begin{eqnarray}\label{kimkimeq1}
div\left(f\nabla f\right)+\lambda f^{2}\leq div\left(f\nabla f\right)+(m-2)|\nabla f|^{2}+\lambda f^{2}=0.
\end{eqnarray}
Considering $M$ a compact Riemannian manifold, integrating (\ref{kimkimeq1}) we have
\begin{eqnarray}\label{kimkimeq2}
\int_{M}\lambda f^{2}dv=\int_{M}\left(div\left(f\nabla f\right)+\lambda f^{2}\right)dv\leq 0.
\end{eqnarray}
Therefore, from (\ref{kimkimeq2}) we can infer that
\begin{eqnarray}\label{kimkimeq3}
\lambda\int_{M}f^{2}dv\leq 0.
\end{eqnarray}
This implies that, either $\lambda\leq0$ or $f$ is trivial. 
It is worth to point out that compact quasi-Einstein metrics on compact manifolds with $\lambda\leq0$ are trivial (see Remark 6 in \cite{kimkim}). 
\hfill $\Box$

\

\noindent {\bf Proof of Theorem \ref{teo2}:}
Let $p$ be a maximum point of $f$ on $M$. Therefore, $f(p)>0$, $(\nabla f)(p)=0$ and $(\Delta f)(p)\geq0$. By hypothesis $R+\lambda(m-n)\leq0$, then from (\ref{oi}) we get
\begin{eqnarray*}
|\nabla f|^{2}\geq\frac{\mu}{m-1}.
\end{eqnarray*}
Whence, in $p\in M$ we obtain
\begin{eqnarray*}
0=|\nabla f|^{2}(p)\geq\frac{\mu}{m-1}.
\end{eqnarray*}
Since $\mu$ is constant, we have that $\mu\leq0$. Moreover, from the third equation in (\ref{test1}) we have
\begin{eqnarray*}
\lambda f^{2}(p)\leq (f\Delta f)(p)+(m-1)|\nabla f|^{2}(p)+\lambda f^{2}(p)=\mu\leq0.
\end{eqnarray*}
Implying that $\lambda\leq0$. Then, from \cite{kimkim} the result follows.

Now, if $R+\lambda(m-n)=0$ from (\ref{oi}) we have that
\begin{eqnarray*}
|\nabla f|^{2}=\frac{\mu}{m-1}.
\end{eqnarray*}
Then, for $p\in M$ we obtain
\begin{eqnarray*}
0=|\nabla f|^{2}(p)=\frac{\mu}{m-1}.
\end{eqnarray*}
Therefore, since $\mu$ is a constant we get that $\mu=0$, i.e., $N$ is Ricci-flat.

\hfill $\Box$

It is worth to say that if $M$ is a compact Riemannian manifold and the scalar curvature $R$ is constant, then $f$ is trivial (see \cite{Case}).

\

\noindent {\bf Proof of Theorem \ref{teo3b}:}
Considering $\lambda=\frac{R}{n-1}$ in equation (\ref{oi}) we obtain
\begin{eqnarray}\label{ooi}
|\nabla f|^{2}+\frac{R}{m(n-1)}f^{2}=\frac{\mu}{m-1}.
\end{eqnarray}

Then, taking the Laplacian we get
\begin{eqnarray}\label{3b1}
\frac{1}{2}\Delta|\nabla f|^{2}+\frac{R}{m(n-1)}\left(|\nabla f|^{2}+f\Delta f\right)=0.
\end{eqnarray}

Moreover, when we consider that $\lambda=\frac{R}{n-1}$ in (\ref{test1}), and contracting the first equation of the system we have that
\begin{eqnarray}\label{3b2}
-\Delta f=\frac{Rf}{m(n-1)}.
\end{eqnarray}
From (\ref{3b2}), (\ref{3b1}) became 
\begin{eqnarray}\label{3b3}
\frac{1}{2}\Delta|\nabla f|^{2}+\frac{R}{m(n-1)}|\nabla f|^{2}=\frac{R^{2}f^{2}}{m^{2}(n-1)^{2}}.
\end{eqnarray}

The first equation of (\ref{test1}) and (\ref{ooi}) allow us to infer that
\begin{eqnarray*}
\frac{2f}{m}Ric(\nabla f)&=&\frac{2Rf}{m(n-1)}\nabla f+2\nabla^{2}f(\nabla f)\nonumber\\
&=&\nabla\left(|\nabla f|^{2}+\frac{Rf^{2}}{m(n-1)}\right)=\nabla\left(\frac{\mu}{m-1}\right)=0.
\end{eqnarray*}
And since $f>0$ we get
\begin{eqnarray}\label{3b4}
Ric(\nabla f, \nabla f)=0.
\end{eqnarray}

Remember the Bochner formula
\begin{eqnarray}\label{bochner}
\frac{1}{2}\Delta|\nabla f|^{2}=|\nabla^{2}f|^{2}+Ric(\nabla f,\nabla f)+g(\nabla f,\nabla\Delta f).
\end{eqnarray}
Whence, from (\ref{3b2}), (\ref{3b4}) and (\ref{bochner}) we obtain
\begin{eqnarray}\label{bochner1}
\frac{1}{2}\Delta|\nabla f|^{2}+\frac{R}{m(n-1)}{|\nabla f|}^{2}=|\nabla^{2}f|^{2}.
\end{eqnarray}
Substituting (\ref{3b3}) in (\ref{bochner1}) we get
\begin{eqnarray}\label{hessiannorm}
|\nabla^{2}f|^{2}=\frac{R^{2}f^{2}}{m^{2}(n-1)^{2}}.
\end{eqnarray}
From the first equation of (\ref{test1}), a straightforward computation give us
\begin{eqnarray}\label{ricnorm}
|Ric|^{2}=\frac{m^{2}}{f^{2}}|\nabla^{2}f|^{2}+\frac{2mR\Delta f}{(n-1)f}+\frac{nR^{2}}{(n-1)^{2}}.
\end{eqnarray}
Finally, from (\ref{hessiannorm}), (\ref{3b2}) and (\ref{ricnorm}) we have that
\begin{eqnarray*}
|Ric|^{2}=\frac{R^{2}}{n-1}.
\end{eqnarray*}
Then, we get the result.
\hfill $\Box$

\

In what follows, we consider the conformal structure given in Section \ref{CFSI} to prove Theorem \ref{teo3a} and Theorem \ref{teo3}.

\

\noindent {\bf Proof of Theorem \ref{teo3a}:}
From definiton, 
\begin{eqnarray}\label{grad}
|\bar{\nabla}f|^{2}=\displaystyle\sum_{i,j}\varphi^{2}\varepsilon_{i}\delta_{ij}f_{x_{i}}f_{x_{j}}=\left(\displaystyle\sum_{i}\varepsilon_{i}\alpha_{i}^{2}\right)\varphi^{2}(f')^{2}=\kappa\varphi^{2}(f')^{2},
\end{eqnarray}
where $\bar{\nabla}f$ is the gradient of $f$ for $\bar{g}$, and $\kappa\neq0$.
Then, from (\ref{eqtop}) and (\ref{grad}) we have
\begin{eqnarray}\label{sera}
\kappa\varphi^{2}(f')^{2}+\left[\frac{\lambda(m-n)+\bar{R}}{m(m-1)}\right]f^{2}=\frac{\mu}{m-1}.
\end{eqnarray}

Consider that $N$ is a Ricci-flat semi-Riemannian manifold, i.e., $\mu=0$, from (\ref{sera}) we get
\begin{eqnarray}\label{sera1}
\frac{f'}{f}=\frac{G(\bar{R})}{\varphi},
\end{eqnarray}
where $G(\xi)=\pm\sqrt{\frac{\kappa[\lambda(n-m)-\bar{R}]}{m(m-1)}}$.
Which give us (\ref{sera3}).

Now, from (\ref{sera1}) we have
\begin{eqnarray}\label{sera2}
\frac{f''}{f}=\left(\frac{G}{\varphi}\right)'+\left(\frac{G}{\varphi}\right)^{2}=\left(\frac{G}{\varphi}\right)^{2}+\frac{G'}{\varphi}-\frac{G\varphi'}{\varphi^{2}}.
\end{eqnarray}

Therefore, from (\ref{system}), (\ref{sera1}) and (\ref{sera2}) we get (\ref{system2}).

\hfill $\Box$

\

\noindent {\bf Proof of Theorem \ref{teo3}:}
Considering that $\bar{R}$ is constant, from (\ref{system2}) we obtain
\begin{eqnarray}\label{system12}
\left\{
\begin{array}{lcc}
(n-2)\varphi\varphi''-mG\varphi'=mG^{2}\\\\
\varphi\varphi''-(n-1)(\varphi')^{2}+mG\varphi'=\kappa\lambda \\\\
(n-1)G\varphi'-mG^{2}=\kappa\lambda
\end{array}
.\right.
\end{eqnarray}
The third equation in (\ref{system12}) give us that $\varphi$ is an affine function. Moreover, since
\begin{eqnarray}\label{hum}
\varphi'(\xi)=\frac{\kappa\lambda+mG^{2}}{(n-1)G},
\end{eqnarray}
we get $\varphi''=0$. Then, from the first and second equations in (\ref{system12}) we have, respectively,
\begin{eqnarray*}
-mG\varphi'=mG^{2}\quad\mbox{and}\quad -(n-1)(\varphi')^{2}+mG\varphi'=\kappa\lambda.
\end{eqnarray*}
This implies that
\begin{eqnarray}\label{hdois}
-(\varphi')^{2}=\frac{\kappa\lambda+mG^{2}}{(n-1)}.
\end{eqnarray}
Then, from (\ref{hum}) and (\ref{hdois}) we get
\begin{eqnarray*}
(\varphi')^{2}+G\varphi'=0.
\end{eqnarray*}
That is, $\varphi'=0$ or $\varphi'=-G$. 

First consider that $\varphi'=0$. From (\ref{scalarcurvature}) and (\ref{system12}), it is easy to see that $\lambda=\bar{R}=0$. Then, we get the first item of the theorem since, as mentioned, the case $\varphi' = 0$ was not considered in \cite{Romildo}.

Now, we take $\varphi'=-G$. Integrating over $\xi$ we have
\begin{eqnarray}\label{phii}
\varphi(\xi)=-G\xi+C;\quad\mbox{where}\quad G\neq0, C\in\mathbb{R}.
\end{eqnarray}
Then, from (\ref{hum}) we obtain
\begin{eqnarray}\label{htres}
\frac{\kappa\lambda+mG^{2}}{(n-1)G}=-G.
\end{eqnarray}
Since $G^{2}=\frac{\kappa[\lambda(n-m)-\bar{R}]}{m(m-1)}$, from (\ref{htres}) we obtain
\begin{eqnarray}\label{scalarcurvature1}
\bar{R}=\frac{n(n-1)\lambda}{(m+n-1)}.
\end{eqnarray}
Considering that $\lambda\neq0$, we can see that $\bar{R}$ is a non-null constant. On the other hand, since $\varphi'=-G$, from (\ref{scalarcurvature}) we get
\begin{eqnarray}\label{anem}
\bar{R}=-n(n-1)\kappa G^{2},
\end{eqnarray}
where $G^{2}=\frac{\kappa[\lambda(n-m)-\bar{R}]}{m(m-1)}$. Observe that (\ref{scalarcurvature1}) and (\ref{anem}) are equivalent. 

Furthermore, from (\ref{sera3}) and (\ref{phii}) we get
\begin{eqnarray*}
f(\xi)=\frac{\Theta}{-G\xi+C}.
\end{eqnarray*}
Now the demonstration is complete.
\hfill $\Box$

\

\noindent {\bf Proof of Corollary \ref{coro1}:}
It is a direct consequence of Theorem \ref{teo3}-(1).

\hfill $\Box$

\

\noindent {\bf Proof of Corollary \ref{coro2}:}
Remember that $\xi=\displaystyle\sum_{i}\alpha_{i}x_{i}$, where $\alpha_{i}\in\mathbb{R}$ (cf. Section \ref{CFSI}). Consider in Theorem \ref{teo3}-(2) that $\alpha_{n}=\frac{1}{G}$ and $\alpha_{i}=0$ for all $i\neq n$. Moreover, taking $C=0$ we get 
\begin{eqnarray}
f(\xi)=\frac{1}{x_{n}^{2}}.
\end{eqnarray}
Moreover, take $\mathbb{R}^{n^{\ast}}_{+}=\{(x_{1},\ldots,x_{n})\in\mathbb{R}^{n}; x_{n}>0\}$. Then, $\left(\mathbb{R}^{n^{\ast}}_{+},g_{can}=\frac{\delta_{ij}}{x_{n}^{2}}\right)=(\mathbb{H}^{n},g_{can})$ is the hyperbolic space. We pointed out that $\mathbb{H}^{n}$ with this metric has constant sectional curvature equal to $-1$. Then, from (\ref{scalarcurvature1}) we obtain $\lambda=-\frac{m+n-1}{n(n-1)}$, and the result follows.

\hfill $\Box$

\

\begin{acknowledgement}
The authors would like to express their deep thanks to professor Ernani Ribeiro Jr for valuable
suggestions.
\end{acknowledgement}

\end{document}